\documentclass[reqno]{amsart}

\usepackage{amsmath}
\usepackage{amssymb}
\usepackage[mathscr]{eucal}
\usepackage{amsthm}
\usepackage{verbatim}
\usepackage[hyphens]{url}
\usepackage{natbib}

\theoremstyle{plain}
\newtheorem{theorem}{Theorem}[section]
\newtheorem{proposition}[theorem]{Proposition}
\newtheorem{lemma}[theorem]{Lemma}
\newtheorem{corollary}[theorem]{Corollary}

\theoremstyle{definition}

\theoremstyle{remark}

\DeclareMathOperator{\Rng}{Rng}

\DeclareMathOperator{\Mix}{Mix}
\DeclareMathOperator{\id}{id}

\newcommand{\mfM}{\mathfrak{M}}

\newcommand{\mcF}{\mathcal{F}}

\newcommand{\mcM}{\mathcal{M}}

\newcommand{\mbR}{\mathbb{R}}

\newcommand{\impl}{~\Rightarrow~}

\begin{document} 

\title{Continuity of the Mixing Operator}

\author{Mikhail Kovtun}
\address{Mikhail Kovtun: Duke University, CDS\\
         2117 Campus Drive, Durham, NC 27708}
\email{mkovtun@cds.duke.edu}
\thanks{This work was supported by
        the Office of Vice Provost for Research, Duke University}


\begin{abstract} 
Mixed distributions are considered as a results of application
of a linear operator, which maps mixing measures to mixed measures.
The main result is a proof of continuity of this mixing operator.
Corollaries for parametric families of distributions
(usually considered in literature) are also discussed.
\end{abstract} 

\maketitle 

\section{Introduction} 
\label{sec:Introduction}

Mixed distributions appear for the first time near the end of 19th
century, in the works of
\cite{Newcomb:1886} and \cite{Pearson:1894}.
The first problem solved in this context was
finding a mixture of two normal distributions.
In the middle of the 20th century, interest in mixed distributions
was renewed in works of \cite{Robbins:1948}
and \cite{Teicher:1960, Teicher:1961, Teicher:1963},
where basic definitions were made and fundamental properties
proved.

These and subsequent works were concentrated on mixtures of
parametric families of distributions.
In the parametric approach, a mixing distribution is defined on
the parameter space.
Alternatively, the mixing measure may be considered
as defined on the space of measures being mixed.
This allows us to consider operation of mixing
as an operator $\Mix : \mfM(\mfM(\Omega)) \rightarrow \mfM(\Omega)$,
where $\mfM(\Omega)$ is the space of measures being mixed
(and where mixtures take values), and $\mfM(\mfM(\Omega))$
is the space of mixing measures.
This operator is obviously linear;
the proof of its continuity is the main topic of
the present article.

To proceed in this way, one needs to specify a class of
basic spaces $\Omega$, a class of $\sigma$-algebras,
and a class of measures.
We propose to consider Radon measures on Borel $\sigma$-algebras
on polish spaces.
This choice, on the one hand, covers most of practically important
cases, and, on the other hand, provides a natural way to obtain
all required structures: the topology
on basic space defines a narrow topology
on the first-level space of measures, which, in turn,
provides a Borel $\sigma$-algebra for the second-level
measures and defines a topology on the second-level space.

We introduce required constructions and list related facts in
section \ref{sec:SpacesOfMeasures}. A proof of continuity
of operator $\Mix$ is given in section \ref{sec:ContinuityOfOperator}.
Section \ref{sec:ParametricFamilies} discusses corollaries
for parametric families.


\section{Spaces of measures} 
\label{sec:SpacesOfMeasures}

For a topological space $\Omega$, $C(\Omega)$ denotes
the space of continuous real-valued functions on $\Omega$,
and $C^*(\Omega)$ denotes its dual, a space of continuous
real-valued linear functionals on $C(\Omega)$.

If $\Omega$ is a Hausdorff topological space,
the space of Radon charges (signed measures) $\mfM(\Omega)$
may be imbedded into $C^*(\Omega)$.
We are interested in the narrow topology on $\mfM(\Omega)$,
i.e. the weakest topology in which functionals
$\phi_f : \mu \mapsto \int f(\omega) \,\mu(d\omega)$ are
continuous for all $f \in C(\Omega)$.
(This topology is usually called ``weak'' in probability theory;
we use term ``narrow'' to distinguish it from the weak topology
on $C^*(\Omega)$.)

When $\Omega$ is a compact metric space, the following strong
facts take place:

\begin{enumerate}
\item
    $C(\Omega)$ is separable metric space in norm topology
    \citep[][IV.13.16]{Dunford:1958}.
\item
    $\mfM(\Omega)$ is isomorphic to $C^*(\Omega)$
    \citep[Riesz representation theorem;][IV.6.3]{Dunford:1958},
    and the narrow topology on $\mfM(\Omega)$ coincides with the
    weak* topology on $C^*(\Omega)$.
\item
    The weak* topology on $C^*(\Omega)$ is metrizable on every
    bounded in norm subset of $C^*(\Omega)$
    \citep[][IV.3.4.4]{Kolmogorov:1972}.
\item
    Every bounded and closed in norm topology subset of $C^*(\Omega)$
    is compact in weak* topology
    \citep[][IV.3.4.5]{Kolmogorov:1972}.
\item
    There exists a countable subset of $C^*(\Omega)$ possessing
    the following property: its intersection
    with every ball $B_r$ is everywhere dense in $B_r$ in weak* topology
    (follows from the compactness of the balls $B_r$).
\end{enumerate}

However, the requirement for $\Omega$ to be compact is too restrictive,
as many interesting spaces (from applicational point of view)
are excluded from the consideration. One possibility to relax
this restriction is to allow $\Omega$ to be a polish
(i.e., complete separable metric) space.

In the case of polish space, neither $C(\Omega)$ is separable,
nor does $C^*(\Omega)$ possess good properties. However,
sufficiently good properties of $\mfM(\Omega)$ can be established
by the following considerations \citep[][chapter III]{Dellacherie:1978}.

As a separable metric space, $\Omega$ has a metric compactification
$b\Omega$ \citep[][2.22.II]{Kuratowski:1966}.
Thus, $\mfM(b\Omega)$ possesses all the above properties.
The natural embedding of $\mfM(\Omega)$ into $\mfM(b\Omega)$
(a charge on $\Omega$ corresponds to a charge on $b\Omega$
carried by $\Omega$) has an important property:
the topology on $\mfM(\Omega)$ generated by $C(\Omega)$ coincides
with the narrow topology inherited from $\mfM(b\Omega)$
\citep[][theorem III.58]{Dellacherie:1978}.
Every bounded closed subset of $\mfM(\Omega)$ is also complete
\citep[][theorem III.60]{Dellacherie:1978}.

Thus, every bounded in norm closed subset $\mcM \subseteq \mfM(\Omega)$
is a polish space in the narrow topology
(but not necessarily compact metric space).
In particular, all probabilistic measures
compose a polish space.
This allows us to consider a space of Radon charges on $\mcM$,
i.e., $\mfM(\mcM)$; it will possess the same topological
properties as $\mfM(\Omega)$.
For our purposes, $\mfM(\Omega)$ serves as a first-level space
of measures, which contains measures being mixed,
and $\mfM(\mcM)$ serves as a second-level space, which
contains mixing measures.

Below, $I_A$ denotes the indicator function of the set $A$,
$\partial A$ denotes a topological boundary of the set $A$,
and $\mfM^+(\Omega)$ denotes a space of non-negative measures on
$\Omega$.

The mixing operator $\Mix : \mfM(\mcM) \rightarrow \mfM(\Omega)$
is defined as:

\begin{equation}   
\label{eq:MixDef}
\Mix(\nu)(A) = \int_\mcM \mu(A) \,\nu(d\mu) \qquad
\text{for every Borel set } A \subseteq \Omega
\end{equation}   

To ensure the correctness of this definition, one needs to
prove the existence of the integral on the right-hand side.
The latest follows from the fact that the mapping
$\mu \mapsto \int f(\omega) \,\mu(d\omega)$ is Borel
for every Borel function $f$ \citep[][III.55,60]{Dellacherie:1978}:
one has $\mu(A) = \int_\Omega I_A(\omega) \,\mu(d\omega)$;
thus, $\mu(A)$ as a function of $\mu$ is Borel, which ensures
existense of the ingegral in (\ref{eq:MixDef}).



\section{Continuity of operator $\Mix$} 
\label{sec:ContinuityOfOperator}

For the sake of simplicity, we restrict ourselves to the
case $\mcM \subseteq \mfM^+(\Omega)$, and prove continuity
of $\Mix$ only on positive cone $\mfM^+(\mcM)$.
Such a restriction, on the one hand, eliminates the necessity
of technical details, while, on the other hand, still covers
the most important case of probabilistic mixtures of probabilistic
measures.

Our proof of continuity of operator $\Mix$ is based on the
following lemma \citep[][III.57]{Dellacherie:1978}:

\begin{lemma}
\label{lm:DM}
Let $f : \Omega \rightarrow \mbR$ be a bounded Borel function.
If a charge $\lambda$ is carried by the set of points of continuity
of $f$, the mapping $\mu \mapsto \int f(\omega) \,\mu(d\omega)$
is continuous at the point $\lambda$.
\end{lemma}

We also need:

\begin{corollary}
\label{cr:DM}
For every Borel $A \subseteq \Omega$,
if a measure $\mu'$ is a point of discontinuity of mapping 
$\phi_A : \mu \mapsto \mu(A)$, then $\mu'(\partial A) \neq 0$. 
\end{corollary}

\begin{proof}
One has $\phi_A(\mu) = \int I_A(\omega) \,\mu(d\omega)$,
and $\partial A$ is the set of points of discontinuity of $I_A$.
If, in contrary, $\mu'(\partial A) = 0$, then by lemma \ref{lm:DM}
$\mu'$ is a point of continuity of $\phi_A$, which contradicts the
assumption.
\end{proof}

\begin{theorem}
\label{th:Main}
For every bounded in norm closed set $\mcM \subseteq \mfM^+(\Omega)$,
operator $\Mix : \mfM^+(\mcM) \rightarrow \mfM^+(\Omega)$ is continuous.
\end{theorem}

\begin{proof}
As both $\mfM(\Omega)$ and $\mfM(\mcM)$ are separable metric spaces,
it is sufficient to show that for every narrowly converging sequence
$\{\nu_k\}_k$ in $\mfM(\mcM)$, $\nu_k \rightarrow \nu_0$,
its image converges to the image of its limit,
$\Mix(\nu_k) \rightarrow \Mix(\nu_0)$.

For this, it is necessary and sufficient to show that for every
Borel $A \subseteq \Omega$ satisfying
\begin{equation}   
\label{eq:bCond}
\Mix(\nu_0)(\partial A) = 0
\end{equation}   

\noindent
one has 
\begin{equation}   
\label{eq:AConv}
\Mix(\nu_k)(A) \rightarrow \Mix(\nu_0)(A)
\end{equation}   

The convergence (\ref{eq:AConv}) means that the mapping
$\nu \mapsto \int_\mcM \mu(A) \,\nu(d\mu)$ is continuous
at the point $\nu_0$. By lemma \ref{lm:DM}, one needs to show
that $\nu_0$ is carried by the set of points of continuity
of mapping $\phi_A : \mu \mapsto \mu(A)$.

Let $M_A$ be the set of points of discontinuity of mapping $\phi_A$.
Suppose, in contrary, that $\nu_0(M_A)>0$.
By corollary \ref{cr:DM}, $\mu \in M_A \impl \mu(\partial A) > 0$;
thus, $\int_\mcM \mu(\partial A) \,\nu_0(d\mu) > 0$.
But $\int_\mcM \mu(\partial A) \,\nu_0(d\mu) = \Mix(\nu_0)(\partial A)$,
and we come to contradiction with condition (\ref{eq:bCond}).
Thus, the convergence (\ref{eq:AConv}) takes place,
and the theorem is proved.
\end{proof}


\section{Parametric families} 
\label{sec:ParametricFamilies}

In parametric approach, one considers a measurable space of parameters
$(\Theta,\mcF)$ and a mapping $\psi : \Theta \rightarrow \mfM(\Omega)$;
the required property is that for every Borel $A \subseteq \Omega$:

\begin{equation}   
\label{eq:mProp}
\psi_A ~:~ \Theta \rightarrow \mbR ~:~ \theta \mapsto \psi(\theta)(A)
\qquad \text{is measurable}
\end{equation}   

The approach considered above can be reduced to the parametric one
by taking $\Theta = \mfM^+(\mcM)$ and $\psi = \id_{\mfM^+(\mcM)}$.
The condition (\ref{eq:mProp}) follows from the fact that
the mapping $\mu \mapsto \mu(A) = \int I_A(\omega) \,\mu(d\omega)$
is Borel for every Borel $A$ (see end of the section 
\ref{sec:SpacesOfMeasures}).

In turn, parametric approach can be reduced to the one considered above
by taking $\mcM = \Rng(\psi)$. The following fact is a straightforward
corollary of requirement (\ref{eq:mProp}):

\begin{proposition}
\label{pr:mPhi}
The mapping $\psi$ is measurable w.r.t. Borel $\sigma$-algebra
of $\mcM$.
\end{proposition}

In many practical cases, mapping $\psi$ is even continuous
(for example, parameterization of normal distributions by means of
mean and variance defines a continuous mapping from real half-plane
to the space of probabilistics measures on real line).

Proposition \ref{pr:mPhi} allows us to defined the mapping
$\hat{\psi} : \mfM(\Theta) \rightarrow \mfM(\mcM)$ by letting
$\hat{\psi}(\nu)(A) = \nu(\psi^{-1}(A))$.
When $\Theta$ is a polish space, a narrow topology can be defined
on $\mfM(\Theta)$.
If the mapping $\psi$ is continuous,
the mapping $\hat{\psi}$ is continuous, too
\citep[][theorem III.8.2]{Shiryaev:2004}.
Thus, theorem \ref{th:Main} implies:

\begin{theorem}
The operator $\Mix_\Theta = \Mix \circ \hat{\psi} : 
\mfM^+(\Theta) \rightarrow \mfM^+(\Omega)$ is continuous,
whenever the mapping $\psi$ is continuous.
\end{theorem}

\bibliographystyle{biostatistics}
\bibliography{math,Probability,Mixtures}

\end{document}